\documentclass{ansarticle}

\usepackage{amsfonts}
\usepackage{amsmath}
\usepackage{enumerate}
\usepackage{graphicx}
\usepackage{color}
\usepackage{array}
\usepackage[english]{babel}
\usepackage{float}

\makeindex

\usepackage{listings}
\definecolor{listinggray}{gray}{0.9}
\definecolor{darkgreen}{rgb}{0,0.7,0}
\definecolor{lbcolor}{rgb}{0.9,0.9,0.9}
\usepackage{amsmath,amssymb}
\usepackage{amsfonts}
\usepackage{theorem}

\usepackage{graphicx}
\usepackage{color}

\usepackage{enumerate}
\usepackage[mathscr]{eucal}
\usepackage{mathrsfs}

\newcommand{\Grad} {\boldsymbol{\nabla}}

\newcommand{\uu}[1]{\boldsymbol{#1}}
\newcommand{\U}{{\uu{u}}}

\renewcommand{\phi}{\varphi}

\title{The LifeV library: engineering mathematics beyond the proof of concept}
\runningtitle{The LifeV library}
\runningauthor{Bertagna, Deparis, Formaggia, Forti, Veneziani}
\author[1]{Luca Bertagna}
\affil[1]{Department of Mathematics and Computer Science, Emory University, Atlanta (GA) 30322 USA}
\author[2]{Simone Deparis\thanks{Corresponding author, \tt{simone.deparis@epfl.ch}}}
\affil[2]{CMCS--MATHICSE--SB, \'Ecole Polytechnique F\'ed\'erale de Lausanne, Station 8, Lausanne, CH--1015, Switzerland}
\author[3]{Luca Formaggia}
\affil[3]{MOX, Dipartimento di Matematica, Politecnico di Milano, Italy}
\author[2]{Davide Forti}
\author[1]{Alessandro Veneziani}

\begin{document}

\maketitle

\begin{abstract}

LifeV is a library for the finite element (FE) solution of partial differential equations
 in one, two, and three dimensions.
It is  written in C++ and designed to run on diverse parallel architectures, including
cloud and high performance computing facilities. In spite of its academic research nature, meaning a library for the development
and testing of new methods,
one distinguishing feature of  LifeV
is its use on real world problems and it is intended to provide a tool
for many engineering applications.
It has been actually used in computational hemodynamics, including cardiac mechanics and
fluid-structure interaction problems,
in porous media, ice sheets dynamics for both forward and inverse problems.
In this paper we give a short overview of the features of LifeV and its coding paradigms
on simple problems.
The main focus is on the parallel environment which is mainly driven by domain decomposition methods and based on
 external libraries such as MPI, the Trilinos project, HDF5 and ParMetis.

\noindent \emph{Dedicated to the memory of Fausto Saleri.}
\end{abstract}






\section{Introduction}
\label{sec:introduction}

LifeV\footnote{Pronounced ``life five'', the name stands for Library 
for Finite Elements, 5th edition as V is the Roman notation for the 
number 5} is a parallel library written in C++ for the approximation 
of Partial Differential Equations (PDEs) by the finite element method
in one, two and three dimensions.
The project started in 1999/2000 as a collaboration between
the modeling and scientific computing group (CMCS) at EPFL Lausanne
and the MOX laboratory at Politecnico di Milano. Later the
REO and ESTIME groups at INRIA joined the project.
In 2006 the library has been progressively parallelized using
MPI with the Trilinos library suite as back-end interface.
In 2008 the Scientific Computing group at Emory University
joined the LifeV consortium.
Since then, the number of active developers has fluctuated between 20 
and 40 people, mainly PhD students and researchers from the 
laboratories within the LifeV consortium.

LifeV is open source and currently distributed under the LGPL license on github\footnote{\url{https://github.com/lifev}, install instructions at \url{https://bitbucket.org/lifev-dev/lifev-release/wiki/lifev-ubuntu}}, 
and migration to BSD License is currently under consideration.
The developers page is hosted by a Redmine system at
\url{http://www.lifev.org}.
LifeV has two specific aims: (i) it provides tools for developing and testing
novel numerical methods for single and multi-physics problems, and
(ii) it provides a platform for simulations of engineering and, more generally, real 
world problems.
In addition to ``basic'' Finite Elements tools, LifeV also provides data structures and algorithms tailored for specific applications in a variety of fields, including fluid and structure dynamics,
heat transfer, and transport in porous media, to mention a few.
It has already been used in medical and industrial contexts, particularly for cardiovascular simulations, including fluid mechanics, geometrical multiscale modeling of the vascular system, cardiac electro-mechanics and its coupling with the blood flow.
When in 2006 we decided to introduce parallelism, the choice has turned towards available open-source tools:
MPI (mpich or openmpi implementations), ParMETIS, and the Epetra, AztecOO, IFPACK, ML, Belos, and Zoltan packages distributed within Trilinos \cite{trilinos}.

LifeV has benefited from the contribution of many PhD students, almost all of them working on project financed by public funds. We provide a list of supporting agencies hereafter.

In this review article we explain the parallel design 
of the library and provide two examples of how to solve PDEs using
LifeV.
Section~\ref{sec:parallel} is devoted to a description of
how parallelism is handled in the library while in 
Section~\ref{sec:FeaturesAndParad} we discuss the distinguishing 
features and coding paradigms of the library. In Section~\ref{sec:life-poisson} we illustrate how to use LifeV to approximate PDEs by the finite element method, using a simple Poisson problem as an example.
In Section~\ref{sec:Oseen} we show how to approximate unsteady Navier--Stokes equations
and provide convergence, scalability, and timings. We conclude by pointing to some applications of LifeV.

Before we detail technical issues, let us briefly address the natural question when approaching this 
software, namely {\it yet another finite element library}?

No question that the research and the commercial arenas offer a huge variety  of finite element libraries
(or, in general, numerical solvers for partial differential equations)
to meet diverse expectations in different fields of engineering sciences. LifeV is - strictly speaking - no exception. 

Since the beginning, LifeV was intended to address two needs: (i) {\it a permanent playground for new methodologies in computational mechanics};
(ii) {\it a translational tool to shorten the time-to-market of new successful methodologies to real engineering problems}.
Also, over the years, we organized portions of the library to be used for teaching purposes.
It was used as a sort of {\it gray box} tool for instance in Continuing Education initiatives at the Politecnico di Milano,
or in undergraduate courses at Emory University (MATH352: Partial Differential Equations in Action) - see \cite{FormaggiaSaleriVenezianiBook}.

As such, LifeV incorporated since the beginning the most advanced methodological developments on topics of interest for the different groups involved. In fact, state-of-the-art methodologies have been rapidly implemented, particularly in incompressible computational fluid dynamics,
to be tested on problems of real interest, so to quickly assess the real performances of new ideas and their practical impact.
On the other hand, advanced implementation paradigms and efficient parallelization were prioritized, as we will describe in the paper. 

When a code is developed with a strong research orientation, the working force is mainly provided by young and junior scholars on specific projects.
This has required a huge coordination effort, in each group and overall. Stratification of different ideas evolving has sometimes
made the crystallization of portions of the code quite troublesome, also for the diverse background of the developers.
Notwithstanding this, the stimulus of the applications has promoted the development of truly advanced methods
for the solution of specific problems. In particular, the vast majority of the stimuli were provided by computational hemodynamics,
as all the groups involved worked in this field with a strict connection to medical and healthcare institutions.
The result is a library extremely advanced regarding performances and with a sort of unique deep treatment of a specific class of problems.
Beyond the overall approach to the numerical solution of the incompressible Navier-Stokes equations (with either monolithic or algebraic partitioning schemes),
fluid-structure interaction problems, patient-specific image based modeling, defective boundary conditions and geometrical multiscale modeling 
have been implemented and tested in LifeV in an extremely competitive and unique way. The validation and benchmarking
on real applications, as witnessed by several publications out of the traditional field of computational mechanics, makes
it a reliable and efficient tool for modern engineering.
So, it is another finite element library yet with peculiarities that make it a significant - somehow unique - example of modern scientific computing tools


\subsection{Financial support}

LifeV has been supported by the 
6th European Framework Programme (Haemodel Project, 2000-2005, 
PI: A.\ Quarteroni),
the 7th European Framework Programme (VPH2 Project, 2008-2011, ERC 
Advanced Grant  MATHCARD 2009-2014, PI: A.\ Quarteroni),
the Italian MIUR (PRIN2007, PRIN2009 and PRIN12 projects, PI: 
A.\ Quarteroni and L.\
Formaggia), 
the Swiss National Science Foundation
(several projects, from 1999 to 2017, 59230, 61862, 65110, 100637, 
117587, 125444, 122136, 122166,  141034,  PI: A.\ Quarteroni),
the Swiss supercomputing initiatives HP2C and PASC (PI: A.\ Quarteroni),
the Fondazione Politecnico di Milano with Siemens Italy (Aneurisk
Project, 2005-2008, PI: A.\ Veneziani), 
the Brain Aneurysm Foundation
(Excellence in Brain Aneurysm Research, 2010, PI: A.\ Veneziani), 
Emory University Research Committee Projects (2007, 2010, 2015, PI:
A.\ Veneziani), 
ABBOTT Resorb Project (2014-2017, PI: H.\ Samady,
D.\ Giddens, A.\ Veneziani), 
the National Science Foundation (NSF DMS 1419060, 2014-2017, PI:
A.\ Veneziani, NSF DMS 1412963, 2014-2017, PI: A.\ Veneziani, NSF DMS 1620406, 2016-2018, PI: A. \ Veneziani),
iCardioCloud Project (2013-2016, PI: F.\ Auricchio, A.\ Reali,
A.\ Veneziani),
the Lead beneficiary programm (Swiss SNF and German DFG, 140184,
2012-2015, PI: A.\ Klawonn,  A.\ Quarteroni, J.\ Schr\"oder),
and the  French National Research Agency. Some private companies also 
collaborated to both the
support and the development, in particular MOXOFF SpA (2012-2104) exploited the geometric multiscale paradigm for the simulation of an industrial packaging system and and Eni SpA (2011-2014) contributed to the development of the Darcy solver and extended FE capabilities, as well as the use of the library for the simulation the evolution of sedimentary basins~\cite{cervone12:_simul}.

\subsection{Main Contributors}
The initial core of developers was the group of A.\ Quarteroni at MOX, Politecnico
di Milano, Italy and at the Department of Mathematics, EPFL, Lausanne, Switzerland from
an initiative of L.\ Formaggia, J.F.\ Gerbeau, F.\ Saleri and 
A.\ Veneziani.
The group of J.F.\ Gerbau at the INRIA, Rocquencourt, France gave significant contributions from
2000 through 2009 (in particular with M.\ Fernandez and, later, 
M.\ Kern). 
The important contribution of C.\ Prud'homme and G.\ Fourestey during their stays at 
EPFL  and of D.\ Di Pietro at University of Bergamo are acknowledged too.
S.\ Deparis has been the coordinator of the LifeV consortium since 2007.

Here, we limit to summarize the list of main contributors 
who actively developed the library in the last five years. 
We group the names by affiliation.
As some of the authors moved over the years to different institutions, they may be listed with multiple affiliations hereafter.

\noindent
[D.\ Baroli, A.\ Cervone, M.\ Del Pra, N.\ Fadel,
E.\ Faggiano, L.\ Formaggia, A.\ Fumagalli, G.\ Iori, 
R.M.\ Lancellotti,  A.\ Melani, S.\ Palamara, S.\ Pezzuto, 
S.\ Zonca]\footnote{\label{mox}MOX, Politecnico di Milano IT}. 
[L.\ Barbarotta, C.\ Colciago, P.\ Crosetto, 
S.\ Deparis, D.\ Forti, G.\ Fourestey, G.\ Grandperrin, 
T.\ Lassila, C.\ Malossi, 
R.\ Popescu, S.\ Quinodoz, S.\ Rossi, 
R.\ Ruiz Baier, P.\ Tricerri]\footnote{\label{epfl}CMCS, EPFL Lausanne, CH}, M.\ Kern\footnote{INRIA, Rocquencourt, FR},
[S.\ Guzzetti, L.\ Mirabella,
T.\ Passerini, A.\ Veneziani]$^{\text{\footnotesize{\ref{mox},}}}$\footnote{\label{emory}Dept. Math \& CS, Emory Univ, Atlanta GA USA},
U.\ Villa$^{\text{\footnotesize{\ref{emory}}}}$, 
L.\ Bertagna$^{\text{\footnotesize{\ref{emory}, }}}$\footnote{\label{fsu} Department of Scientific Computing, Florida State University, Tallahassee, FL USA},\footnote{\label{snl}Sandia National Lab, Albuquerque, NM USA},
M.\ Perego$^{\text{\footnotesize{\ref{mox},\ref{emory},\ref{fsu},\ref{snl}}}}$,
A. \ Quaini\footnote{Department of Mathematics, University of Houston, TX USA}$^{\text{\footnotesize{,\ref{emory}}}}$,
H.\ Yang$^{\text{\footnotesize{\ref{emory},\ref{fsu}}}}$,
A. \ Lefieux\footnote{Department of Civil Engineering and Structures, University of Pavia, IT}$^{\text{\footnotesize{,\ref{emory}}}}$. 

\paragraph*{Abbreviations}: Algebraic Additive Schwarz (AAS), FE (Finite Elements), FEM (Finite Element Method), DoF (Degree(s) of Freedom), OO (Object Oriented), PCG (Preconditioned Conjugate Gradient), PGMRes (Preconditioned Generalized Minimal Residual), ML (Multi Level), DD (Domain Decomposition), MPI (Message Passing Interface), ET (Expression Template), HDF (Hierarchical Data Format), HPC (High Perfromance Computing), CSR (Compact Sparse Row)


\section{The Parallel Framework}
\label{sec:parallel}

The library can be used for the approximation of PDEs in
one dimension, two dimensions, and three dimensions. Although it can be used in serial mode (i.e., with one
processor), parallelism is crucial when solving three dimensional problems. To better underline the ability of LifeV
to tackle large problems, in this review we focus on PDEs discretized on unstructured linear tetrahedral meshes,
although we point out that LifeV also supports hexahedral meshes as well as quadratic meshes.

Parallelism in LifeV is achieved by domain decomposition (DD) strategies, although it is not
mandatory to use DD preconditioners for the solution of sparse linear systems.
In a typical simulation, the main steps involved in the parallel solution of the finite element problem using LifeV are the following:
\begin{enumerate}
  \item All the MPI processes load the same (not partitioned) mesh.
  \item The mesh is partitioned in parallel using ParMETIS or Zoltan. At the end each process keeps only its own local partition.
  \item The DoFs are distributed according to the mesh partitions. By looping on the local partition, a list of local DoF in global numbering is built.
  \item The FE matrices and vectors are distributed according to the DoFs list. In particular,
  the matrices are stored in row format, for which whole rows are assigned to the process owning the associated DoF.
  \item Each process assembles its local contribution to the matrices and vectors. Successively, global communication consolidates contributions on shared nodes (at the interface of two subdomains).
  \item The linear system is solved using an iterative solver, typically either a Preconditioned Conjugate Gradient (PCG)
  when possible or a Preconditioned GMRes (PGMRes). The preconditioner runs in parallel. Ideally, the number of preconditioned iterations should be independent of the number of processes used.
  \item The solution  is downloaded to mass storage
  in parallel using HDF5 for post-processing purposes (see Sect.~\ref{sec:paraview}).

\end{enumerate}
The aforementioned steps are explained in detail in the next subsections.

\subsection{Mesh partitioning: ParMETIS and Zoltan}
\label{sec:parmetis}

As mentioned above, LifeV achieves parallelism by partitioning the mesh among the available
processes. Typically, this is done ``online'': the entire mesh is loaded by all the processes but
it is deleted after the partitioning phase, so that each process keeps only the part required for the solution of the
local problem and to define inter-process communications. As the mesh size increases, the ``online'' procedure
may become problematic. Therefore for large meshes it is possible, and sometimes
necessary, to partition the mesh
offline on a workstation with sufficient memory~\cite{radu:PhD}. It is also
possible to include an halo of ghost elements such that the partitions overlap by one or more layers of elements (see e.g. \cite{guzz1}). This may be relevant for schemes that require a large stencil. To perform the partition, LifeV can interface with two third party libraries: ParMetis and Zoltan \cite{metis,zoltan}.




\subsection{Distributed arrays: Epetra}
\label{sec:epetra}

The sparse matrix class used in LifeV is a wrapper to the Epetra matrix container
Epetra\_FECrsMatrix and, similarly, the vector
type is a wrapper to Epetra\_FEVector, both provided by the Epetra~\cite{heroux09:_epetra} package of Trilinos.
The distribution of the unknowns is determined automatically by the
partitioned mesh: with a loop over each element of the local mesh we create the list
of DoF managed by the current processor. This  procedure in fact creates a {\it repeated} map,
i.e., an instance of an Epetra\_Map with some entries referring to the DoF associated with geometric entities lying on the  interface between two (or more) subdomains. Then, a {\it unique}
map is  created, in such a way that, among all the owners of a \textit{repeated} DoF, only one will also own it in the \textit{unique} map.
The unique map is used for the vectors and matrices to be used in the linear algebra routines as well as for the solution vector. The repeated map is used
 to access information stored on other processors, which is usually necessary only in the assembly and post-processing phases.


The assembly of the FE matrices is typically performed by looping on the local elements~\cite{ernguermond,FormaggiaSaleriVenezianiBook}.
To reduce latency time, the loops on each subdomain are performed in parallel, without need of any communication
during the loop. Just a single communication phase takes place once all processes have assembled their local contribution, to complete the assembly for interface dofs.

Efficiency and stability may be improved by two further available operations,
(i) precomputing the matrix graph; (ii) using overlapping meshes.
The former demands for the creation at the beginning of the simulation of an Epetra\_Graph, associated to the matrix. Since it depends on the problem at hand and the chosen finite elements,
its computation needs a loop on all the elements. {This is coded by {\it Expression Templates}} (ET), see Sect. \ref{sec:et}, using the same call sequence as for the matrix assembly.
The latter further reduces communications by allowing all processes to compute the local finite element matrix also on all elements sharing a DoF on the interface. As a result, each process can independently compute all the entries of matrices and vectors pertaining to the DoF it owns at the price of some extra computation. Yet, the little overhead is justified by the complete elimination of the post-assembly communication costs.

\subsection{Parallel preconditioners}
\label{sec:prec}
The solution of linear systems in LifeV relies on the Trilinos~\cite{trilinos} packages AztecOO
 and Belos~\cite{bavier12:_amesos_belos}, which provide an extensive choice of iterative or direct
 solvers. LifeV provides a common interface to both of them.

The proper use of ParMETIS and Zoltan for the partitioning, and of Epetra matrices and vectors for the linear algebra, ensures that  matrix-vector multiplication and vector operations are properly parallelized, i.e., they scale well with the number of processes used, and communications are optimized. In this situation the parallel scalability of iterative solvers like PCG or PGMRES depends essentially on the properties of the preconditioner.

The choice of preconditioner is thus critical. In our experience it may follow two directions: (i) 
parallel preconditioners for the generic linear systems, like single or multilevel overlapping Schwarz preconditioners, or multigrid preconditioners, that are generally well suited for highly coercive elliptic problems, or incomplete factorization (ILU), which are generally well suited for advection-dominated elliptic problems; (ii) 
problem specific preconditioners, typically required for multifield or multiphysics problem. These preconditioners exploits specific features of the problem at hand to recast the solution to standard problems that can be eventually solved with the generic strategies in (i). Preconditioners of this class are, e.g., the SIMPLE, the Least Square Commutator, the Caouet-Chabard and the Yosida preconditioners for the incompressible Navier-Stokes equations~\cite{PatankarSpalding1972,elman2014finite,cahouet1988some,Veneziani2003}
 or the Monodomain preconditioner for the Bidomain problem in electrocardiology~\cite{giorda09}.

In LifeV preconditioners for elliptic problems are indeed an interface to the Trilinos package IFPACK~\cite{ifpack}, which is a suite of algebraic preconditioners based on incomplete factorization, and to ML~\cite{ml-guide} or MueLu~\cite{MueLu},
which are two Trilinos packages for multi-level preconditioning based on algebraic multigrid.
Typically we use IFPACK to define algebraic overlapping Schwarz preconditioners with exact or inexact LU factorization
of the restricted matrix. The preconditioner $P$ in this case can be formally written as
\begin{equation}
P^{-1} = \sum_{i=1}^n R_i^T A_i^{-1} R_i, \text{ where }
A_i = R_i A^{-1} R_i^T, \label{eq:schwarz}
\end{equation}
where $A$ is the finite element matrix related to the PDE approximation, $n$ is the number of partitions (or subdomains)
$\Omega_i$, $R_i$ is the restriction operator to $\Omega_i$ and $R_i^T$ the
extension operator from $\Omega_i$ to the whole domain $\Omega$.
$A_i$ is inverted many times during the iterations of PCG or PGMRES, which is why it is factorised
by LU or ILU. In LifeV, the choice of the factorization is left to the user.


Similarly, it is possible to use multilevel preconditioners via the ML~\cite{ml-guide} or MueLu~\cite{MueLu,MueLuURL} packages.
They work at the algebraic level too
and the coarsening and extension are done either automatically, or by user defined strategies, based on the parallel distribution of the matrix and its graph.
The current distribution of LifeV does not offer the last option, but the  interested developer could add this extra functionality with relatively little effort.

\subsection{Parallel I/O with  HDF5}
\label{sec:paraview}

When dealing with large meshes and large number of processes, input/output access to files on disk deserves particular care.
A prerequisite is that the filesystem of the supercomputing architecture being used provides the necessary access speed in parallel. However this is not enough. MPI itself offers parallel I/O capabilities, for which HDF5 is one
of the existing front-ends~\cite{hdf5}. Although other formats are also supported (see section \ref{subsec:io_formats}), LifeV strongly encourages the use HDF5 for I/O processing essentially for three reasons.
\begin{itemize}

\item The number of files produced is independent of the number of running processes: each process accesses the same file in parallel
and writes out his own chunks of data. 
As a result, LifeV generally produces one single large output binary file, along with a xmf text file
describing its contents.

\item HDF5 is compatible with open source post-processing visualization tools like Paraview~\cite{paraview} and  VisIt~\cite{visit}).

\item Having one single binary file makes it very easy to use for restarting a simulation. 
\end{itemize}

The interface to HDF5 in LifeV exploits the facilities of the EpetraExt package in Trilinos, and since the LifeV vectors are compatible with Epetra format the
calls are simple~\cite{heroux09:_epetra}.



\section{Features and Paradigms}\label{sec:FeaturesAndParad}

\subsection{I/O data formats} \label{subsec:io_formats}

In a typical simulation, the user provides a text file containing input data (including physical and discretization parameters, and options for the linear/nonlinear solvers), and the code will generate results, which need to be stored for post-processing. 
Although it is up to the user to write the program main file where the data file is parsed, LifeV makes use of two particular classes in order to forward the problem data to all the objects involved in the simulation: the GetPot class \cite{getpot} (for which LifeV also provides an \emph{ad hoc} re-implementation inside the core module), and the ParameterList class from the Teuchos package in Trilinos. The former has been the preferred way since the early development of LifeV, and is therefore supported by virtually all classes that require a setup. The latter is used mostly for the linear and nonlinear solvers, since it is the standard way to pass configuration parameters to the Trilinos solvers. Both classes map strings representing the names of properties to their actual values (be it a number, a string, or other), and they both allow the user to organize the data in a tree structure. 
Details on the syntax of the two formats
are available online.

When it comes to mesh handling, LifeV has the built-in capability of generating structured meshes on domains of the form $\Omega=[a_1,b_1]\times[a_2,b_2]\times[a_3,b_3]$. If a more general mesh is required, the user needs to create it beforehand. Currently, LifeV supports the FreeFem \cite{freefem} and Gmsh \cite{gmsh} formats (both usually with extension \verb|.msh|) in 2D, while in 3D it supports the formats of Gmsh, Netgen \cite{netgen} (usually with extension \verb|.vol|) and Medit \cite{medit} (usually with extension \verb|.mesh|). Additionally, as mentioned in Section \ref{sec:parmetis}, LifeV offers the capability of offline partitioning.  The partitions of a mesh are stored and subsequently loaded using HDF5 for fast and parallel input.

Finally, LifeV offers three different formats for storing the simulation results for post-processing: Ensight \cite{ensight}, HDF5 \cite{hdf5}, 
which is the preferred format when running in parallel, and VTK \cite{vtk}. All these formats are supported by the most common scientific visualization software packages, like Paraview \cite{paraview} and VisIt \cite{visit}. The details on these formats can be found on their respective webpages.

\subsection{Expression Templates for Finite Elements}
\label{sec:et}

One of the aims of LifeV is to be used in multiple contexts, ranging from industrial and social applications to teaching 
purposes. For this reason, it is important to find the best trade-off between computational efficiency 
and code readability.
The high level of abstraction proper of C++ is in principle perfectly matched by the abstraction of mathematics. 
 However, versatility and efficiency may conflict. Quoting~\cite{furnish},
a ``natural union has historically been more of a stormy relationship". This aspect is crucial in High Performance Computing, where efficiency is a priority.
Operator overloading has a major impact on efficiency, readability, maintainability and versatility, 
however it may adversely affect the run time. 
{\it Expression Templates} (ET) have been originally devised to minimize this drawback  by T. Veldhuizen~\cite{veldhuizen},
and later further developed in the context of linear algebra~\cite{HArdtlein2009,iglberger2012expression} and solution of partial differential equations~\cite{pflaum2001expression,di2009expression}.

In the context of linear algebra, the technique was developed to allow  high level vector syntax without compromising on code speed, due to function overloading. The goal of ET is to write high level \emph{Expressions}, and use \emph{Template} meta-programming to parse the expression at compile-time, generating highly efficient code. Put it simply, in the context of PDE's, ET aims to allow a syntax of the form

\begin{verbatim}
auto weak_formulation = alpha*dot(grad(u), grad(v)) + sigma*u*v;
\end{verbatim}
which is an \emph{expression} very close to the abstract mathematical formulation of the problem. However, during the assembly phase, the resolution of the overloaded operators and functions would yield a performance hit, compared to a corresponding \emph{ad-hoc} for loop. To overcome this issue, the above \emph{expression} is implemented making massive use of \emph{template} meta-programming, which allows to expand the \emph{expression} at compile time, resulting in highly efficient code. Upon expansion, the \emph{expression} takes the form of a combination of polynomial basis functions and their derivative at a quadrature point. At run time, during the assembly phase, such combination is then evaluated \emph{at once} at a given quadrature point, as opposed to the more \emph{classical} implementation, where all contributions are evaluated separately and then summed up together.

For instance, for the classical linear advection-diffusion-reaction equation in the unknown $u$, $-\mu \Delta u + \mathbf{\beta} \cdot \nabla u + \sigma u$,
we need to combine three different differential operators weighted by coefficients $\mu,\mathbf{\beta}$ and $\sigma$.
These may be numbers, prescribed functions or pointwise functions (inherited for instance by another FE computation).
In nonlinear problems --- after proper linearization --- expressions may involve finite element functions too. 
Breaking down the assembly part to each differential operator individually with its own coefficient is possible but leads to  duplicating the loops over
the quadrature nodes, as opposed to the assembly of their sum.
In 
LifeV the possible
differential operators for the construction of linear and nonlinear advection-diffusion-reaction problems are enucleated into
specific \emph{Expression}s,
following the idea originally proposed in \cite{eta}.
The ET technique provides readable code
with no efficiency loss for operator parsing. As a matter of fact, the final gathering of all the assembly operations in a single loop,
as opposed to standard approaches with a separate assembly for each elemental operator, introduces computational advantages. Indeed, numerical tests have pointed out a significantly improved performance for probelms with non-constant coefficients when using the ET technique.

 A detailed description of ET definition and  implementation in LifeV can be found in~\cite{Quinodoz:PhD} and in the code snapshots presented later on. 
\section{Basics: Life = Library of Finite Elements}\label{sec:life-poisson}
The library supports different type of finite elements. The use of ET makes the set-up of simple problems easy, as we illustrate hereafter.

\subsection{The Poisson problem}
\label{sec:PDE}
As a first example, we present the setup of a finite element solver for a Poisson
problem.
We assume a polygonal $\Omega$ with boundary $\partial\Omega$  split into two subsets $\Gamma^D$ and $\Gamma^N$ of positive measure such that
$\Gamma^D\cup\Gamma^N = \partial\Omega$ and $\Gamma^D\cap\Gamma^N=\emptyset$.
Let $V^D_h\subset H^1_{\Gamma^D}$ be a discrete finite element space relative to a mesh
of $\Omega$, for example continuous piecewise linear functions vanishing on $\Gamma_D$.
The Galerkin formulation of the problem reads: find $u_h\in V^D_h$ such that
\begin{equation}
\label{eq:poisson}
\int_{\Omega} \kappa \Grad u_h : \Grad \phi_h
=
- \int_{\Omega} \kappa \Grad u^D  : \Grad \phi_h
+\int_\Omega f \phi_h
+\int_{\Gamma^N} g^N \phi_h
\qquad \forall \phi_h \in V^D_h,
\end{equation}
where $\kappa$ is the diffusion coefficient,
possibly dependent on the space coordinate, $g^N$ is the Neumann boundary condition,
$\frac{\partial u}{\partial n} = g^N$ on $\Gamma^N$, and $f$ are the volumetric forces.
The lifting $u_D$ can be any finite element function such that $u^D|_{\Gamma^D}$
is a suitable approximation of $g^D$. As usually done, $u^D$ is such that $u^D|_{\Gamma^D}$
is the Lagrange interpolation of $g^D$, extended to zero inside $\Omega$.

In LifeV, DoF associated with Dirichlet boundary conditions are not physically eliminated from the FE unknown vectors
and matrices.
Even though this elimination would certainly affect positively the performances of the linear algebra solver,
it introduces a practical burden in the implementation and memory particularly in 3D unstructured problems, that makes it less appealing.
The enforcement of these conditions can be done alternatively in different ways as illustrated e.g. in~\cite{FormaggiaSaleriVenezianiBook},
after the matrix assembly. We illustrate the strategy adopted in the following sections.

\subsection{Matrix form and expression templates}
\label{sec:ETA}
Firstly, we introduce the matrix assembled for homogeneous Natural conditions associated with the
differential operator at hand (aka ``do-nothing'' boundary conditions, as they do not require any extra work to the pure discretization of the differential operator),
i.e.
$$
A =(a_{ij})_{i,j=1,...,n} \text{ and }
a_{ij} =
\int_{\Omega} \kappa \Grad \phi_j : \Grad \phi_i, i,j=1,...,n,
$$
where $\phi_j$, $j=1,...,n$ are the basis functions of the finite element space $V_h$.
Thanks to the ET framework ~\cite{eta,Quinodoz:PhD}, once the mesh, the solution FE space and the quadrature rule have been created, the assembly of the stiffness matrix $A$ is as simple as the following instruction

\begin{lstlisting}[language=C++]
integrate ( elements (localMeshPtr), quadratureRule, uFESpace, uFESpace,
            value(kappa) * dot ( grad ( phi_i ) , grad ( phi_j ) )     )  >> systemMatrixPtr;
\end{lstlisting}

We emphasize how the ET syntax clearly highlights the differential operator being assembled, making the code easy to read and maintain. In a similar way, we define the right hand side of the linear problem as the
vector
$$
\mathbf b = (b_i)_{i=1,...,n}\text{ where }
b_i = \int_\Omega f \phi_i + \int_{\Gamma^N}g^N\phi_i, i=1,...,n.
$$
In this case, possible non-homogeneous Neumann or natural conditions are included.
Finally, the DoF related to the Dirichlet boundary conditions are enforced
by setting the associated rows of $A$ equal to zero except for the diagonal entries.
In this way, the equation associated with the $i$-th Dirichlet DoF is replaced by $c u_i = c g_i$, where  $c$ is a scaling factor, depending in general on the mesh size, to be used to control the condition number of the matrix.
A general strategy is to pick up  values of the same order of magnitude
of the entries of the row of the do-nothing matrix being modified.

Without further modification, the system matrix is not symmetric anymore.
Many of the problem faced in application are not symmetric, therefore we describe here
only how to deal with non-symmetric matrices.

It is worth noting that the symmetry break does not prevent using specific methods for symmetric systems
like CG when appropriate as pointed out in \cite{ernguermond}.
A symmetrization of the matrix can be also  achieved by enforcing the condition $c u_i = c g_i$
column-wise, i.e. by setting to 0 also the off-diagonal entries of the columns of Dirichlet DoF.
Some sparse-matrix formats   oriented to row-wise access of the matrix, like the popular CSR,
need in this case to be equipped with specific storage information that could make the column-wise access convenient~\cite{FormaggiaSaleriVenezianiBook}.

\subsection{Linear algebra}
\label{sec:DD}

The linear system $A\mathbf x = \mathbf b$ can be solved by a preconditioned iterative method like
PCG, PGMres, BiCGStab, etc, available in the packages AztecOO or Belos of Trilinos. The following snippet highlights the simplicity of the usage of LifeV's linear solver interface.

\begin{lstlisting}[language=C++]
// Load the linear solver parameters
Teuchos::RCP<Teuchos::ParameterList> params = Teuchos::getParametersFromXmlFile("params.xml");
// Create and setup the preconditioner
std::shared_ptr<Preconditioner> precPtr (new PreconditionerIfpack);
precPtr->setParametersList (params->sublist("Preconditioner").sublist("Ifpack"));
// Create and setup the linear solver
LinearSolver linearSolver (Comm);         // Comm is the communicator
linearSolver.setOperator (systemMatrix);  // Set the matrix of the system
linearSolver.setParameters (*params);     // Set the linear solver parameters
linearSolver.setPreconditioner (precPtr); // Set the preconditioner
linearSolver.setRightHandSide (rhs);      // Set the right hand side
// Solve the linear system, and store the solution in the given vector
linearSolver.solve (solution);
\end{lstlisting}

The choice of the method and its settings are to be set via an input xml file. For the above example (which used AztecOO as solver and Ifpack as preconditioner), a minimal input file {\tt params.xml} would have the following form

\begin{lstlisting}[language=xml]
<ParameterList>
    <Parameter name="Solver Type" type="string" value="AztecOO"/>
    <Parameter name="Prec Type"   type="string" value="Ifpack"/>
	<ParameterList name="Solver: Operator List">
		<ParameterList name="Trilinos: AztecOO List">
    		<Parameter name="solver"   type="string" value="gmres"/>
	    	<Parameter name="conv"     type="string" value="rhs"/>
    		<Parameter name="tol"      type="double" value="1.e-9"/>
	    	<Parameter name="max_iter" type="int"    value="50"/>
    		<Parameter name="kspace"   type="int"    value="20"/>
    	</ParameterList>
    </ParameterList>
	<ParameterList name="Preconditioner">
		<ParameterList name="Ifpack">
    		<Parameter name="prectype"            type="string" value="ILU"/>
    		<Parameter name="overlap level"       type="int"    value="2"/>
            <Parameter name="fact: level-of-fill" type="int"    value="2"/>
    	</ParameterList>
    </ParameterList>
</ParameterList>
\end{lstlisting}

\begin{figure}[t]
  \centering
  \includegraphics[bb=7 0 537 414,width=0.45\textwidth]{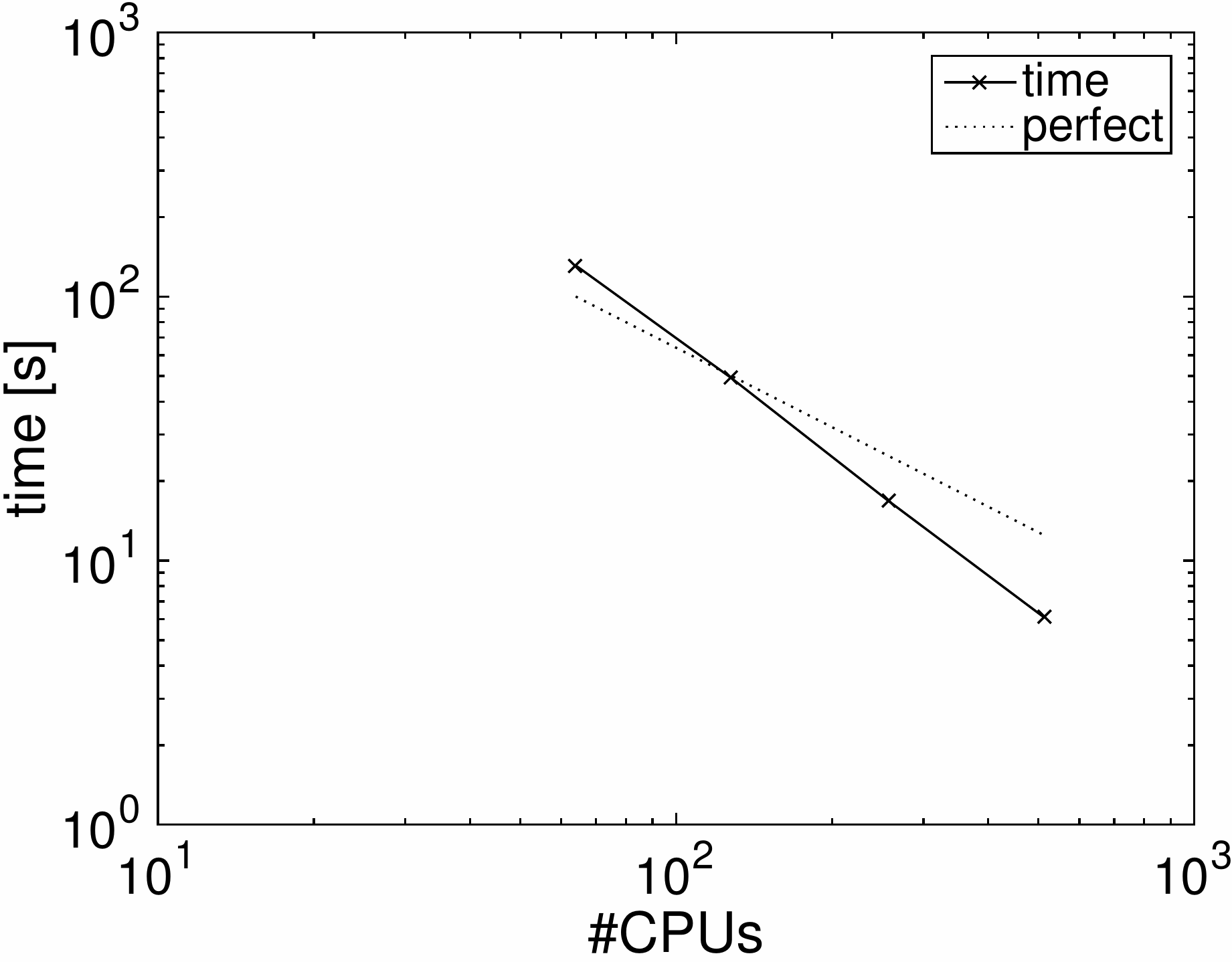}
  \includegraphics[bb=18 0 539 410,width=0.45\textwidth]{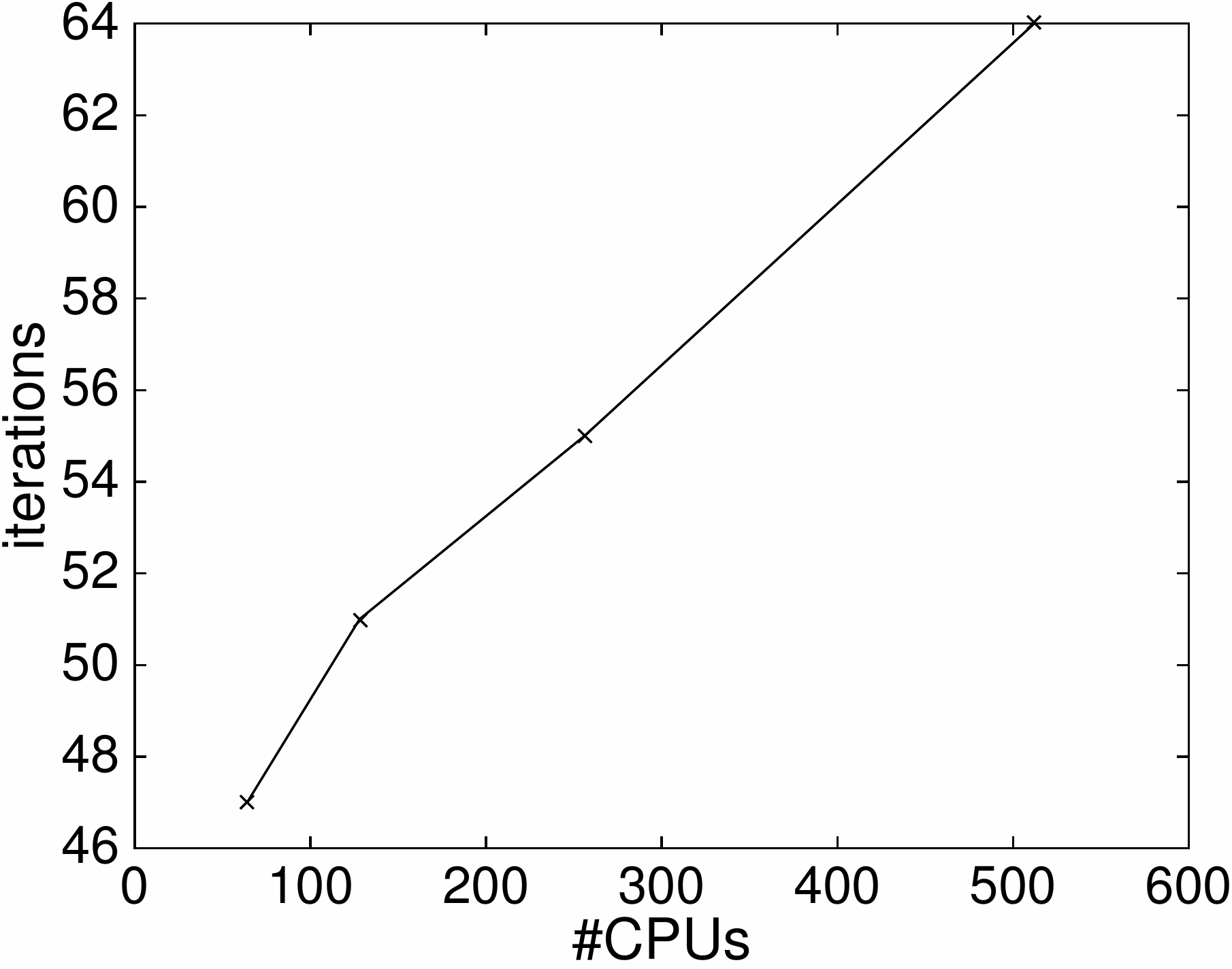}
  \caption{Solving a Poisson problem in a cube with P2 finite elements with
    1'367'631 degrees of freedom. The scalability in terms of CPU time (left) is perfect, however
  the number of iterations (right) linearly increases. The choice of the preconditioner is not optimal, the use
  of a coarse level or of multigrid in the preconditioner is essential and allows to use more processes with no
  loss of resources, cf.\ also Figures~\ref{fig:VMSscalability} and~\ref{fig:FSIscalability}.
}
  \label{fig:laplacian}
\end{figure}

The main difficulty is to set up a scalable preconditioner.  As pointed out , in LifeV there are
several options based on Algebraic Additive Schwarz (AAS) or Multigrid preconditioners. In the first case,
the local problem related to $A_i$ in~\eqref{eq:schwarz} has to be solved.
It is possible to use an LU factorization, using the interface with Amesos~\cite{amesos,amesos:para06} or
incomplete factorizations (ILU).
LU factorizations are more robust than incomplete ones in the sense that they do not need any
parameter tuning, which is delicate in particular within AAS, while incomplete ones are much faster and require less storage.
In a parallel context though, for a given problem,
the size of the local problem is inversely proportional to the number
of subdomains. The LU factorization, whose cost depends only on the number of unknowns,
is perfectly scalable, but more memory demanding. An example of the scalability in this settings is given in Figure~\ref{fig:laplacian}.
LifeV leaves the choice to the user depending on the type of problem and
computer architecture at hand.

The previous example is just an immediate demonstration of LifeV coding. For other examples we refer the reader  to~\cite{FormaggiaSaleriVenezianiBook}.

\section{The CFD Portfolio}
\label{sec:Oseen}

A core application developed since the beginning, consistently with the tradition of the group where the library has been 
originally conceived, is incompressible fluid dynamics, which is particularly relevant 
for hemodynamics.
It is well-known that the problem has a {\it saddle point} nature that stems from the incompressiblity constraint.  
From the mathematical stand point, this introduces specific challenges, for instance the 
choice of finite element spaces for velocity and pressure that should satisfy the so called {\it inf-sup condition}, unless special stabilization 
techniques are used~\cite{elman2014finite}. LifeV offers both possibilities.
In fact, one can choose among inf-sup stable P2-P1 finite element pairs or equal order P1-P1 or P2-P2 stabilized formulations, either
by interior penalty \cite{burmanfernandez} or SUPG-VMS~\cite{Bazilevs2007173,FortiDedeVMS15}.

As it is well known, for high Reynolds flow it is important to be able to describe turbulence by modeling it, being impossible to resolve it in practice. Being hemodynamics the main LifeV application, where turbulence is normally less relevant, LifeV has not implemented a full set of turbulence models. However, it includes the possibility of using the Large Eddie Simulation (LES)
approach, which relies on the introduction of a suitable filter of the
convective field in the Navier-Stokes equations with the role of bringing  the 
unresolved scales of turbulence 
to the mesh scale.
The Van Cittert deconvolution operator considered in~\cite{BowersRebholz}
has been recently introduced in LifeV in~\cite{BertagnaQuainiVeneziani}. 
Validation up to a Reynolds number 6500 has been validated with the FDA Critical Path Initiative Test Case.

Another LES procedure based on the variational splitting of resolved and unresolved 
parts of the solution has been considered in~\cite{FortiDedeVMS15}, while other LES filtering techniques, in particular the $\sigma$-model~\cite{nicoud2011using}, have been implemented in~\cite{Lancellotti15}.

\subsection{Preconditioners for Stokes problem}
\label{preconditioners}

In Section~\ref{sec:prec} we have introduced generic parallel preconditioners based on AAS or multigrid.
These algorithms have been originally devised for elliptic problems. 
In our experience, their use for 
saddle-point problems like Darcy, Stokes and Navier--Stokes equations, or in fluid structure interaction problems, is not effective.
However, their combination with specific preconditioners, like e.g. SIMPLE, Least Square Commutator,  Yosida for unsteady  Navier--Stokes equations
or Dirichlet-Neumann for FSI, leads to efficient and scalable solvers.

\begin{figure}[t]
  \centering
  \includegraphics[bb=0 0 350 330,width=0.45\textwidth]{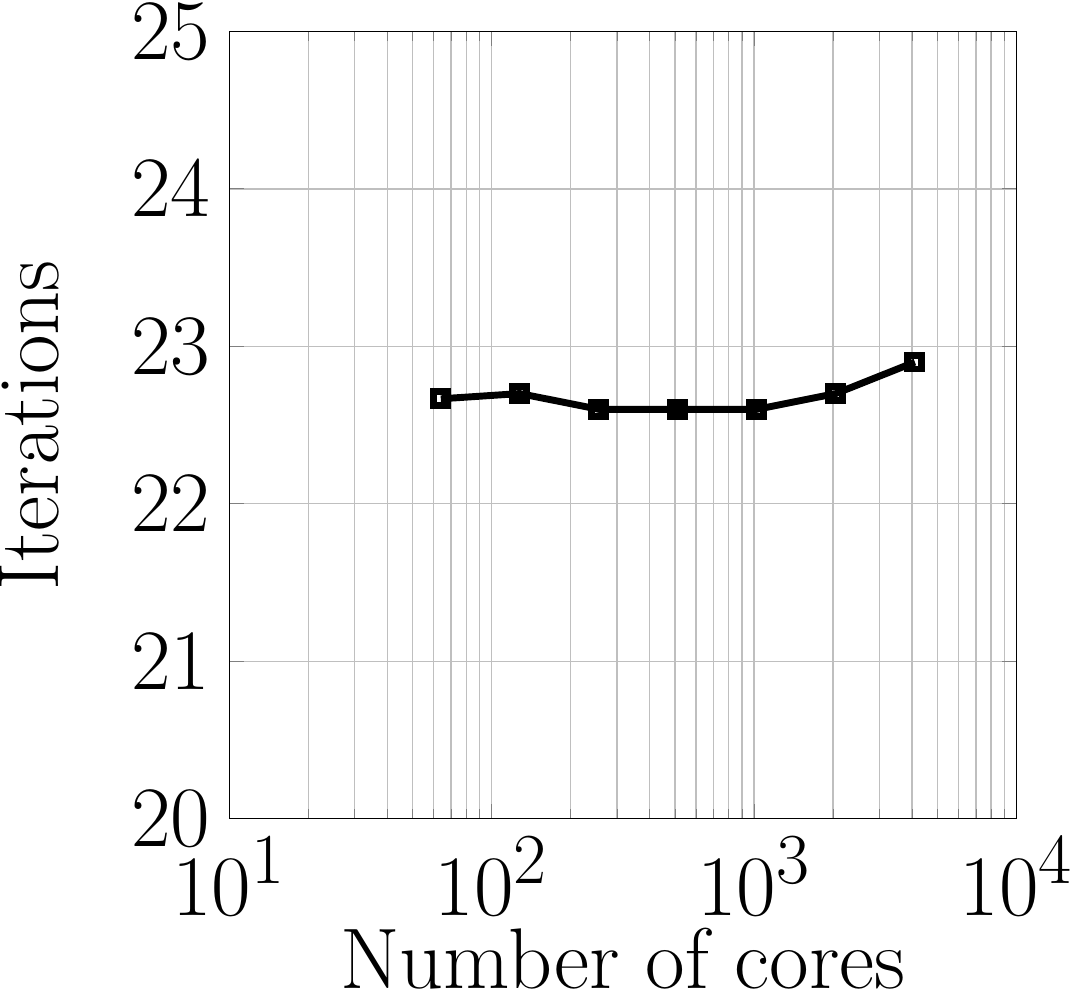} 
  \includegraphics[bb=0 0 350 330,width=0.45\textwidth]{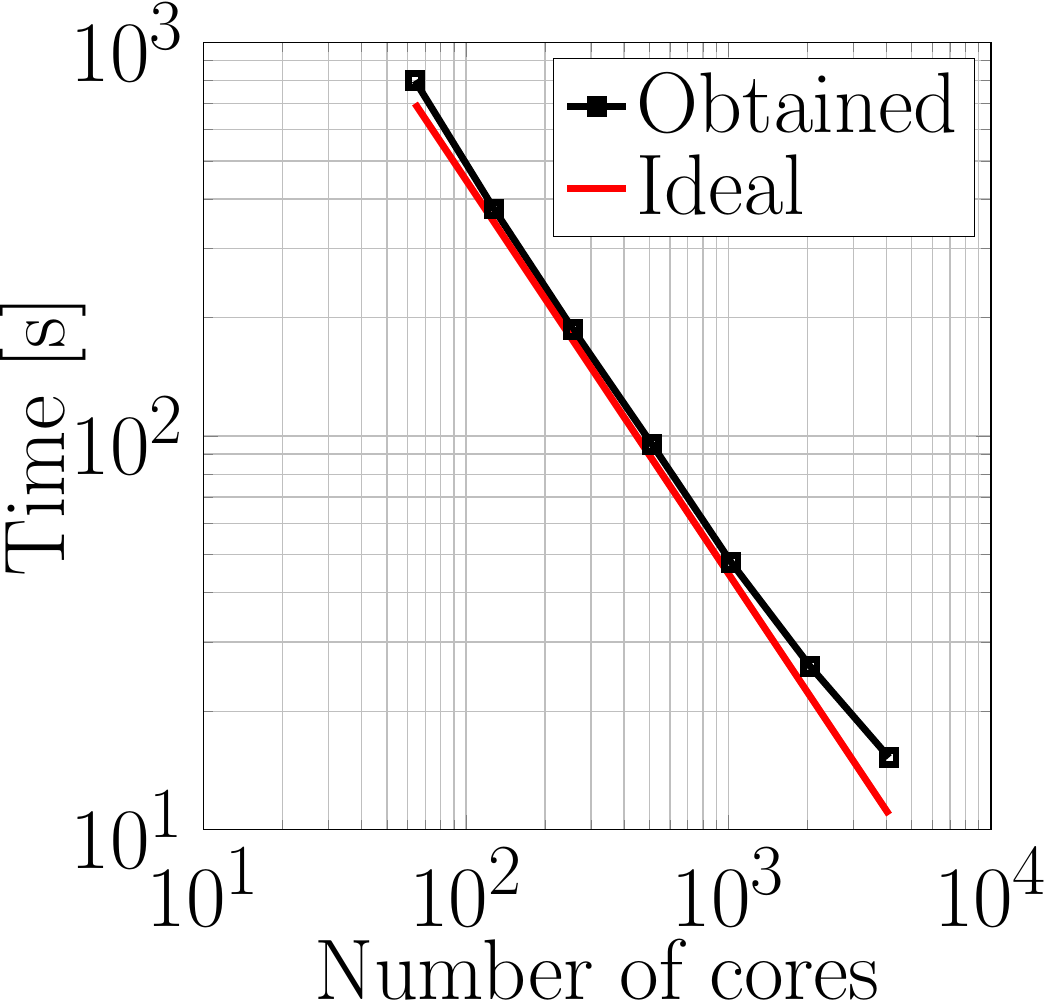}
  \caption{Flow around a cylinder, Reynolds number equal to 22'000.
Scalability of hybrid preconditioners for stabilized Navier--Stokes equations~\cite{FortiDedeVMS15}.
  Simulations on PizDora at CSCS. VMS-LES stabilized P2-P2 finite elements, 9'209'040 degrees of freedom
  and time step 0.0025 s. On the left the number of iterations and on the right the time to solve an entire time step.}
  \label{fig:VMSscalability}
\end{figure}

In~\cite{DeparisGrandperrinQuarteroni2013} this approach has been applied with
success to unsteady Navier-Stokes problems with inf-sup stable finite elements, 
and then extended to a VMS-LES stabilized formulation with equal order elements for 
velocity and pressure~\cite{FortiDedeVMS15}, see also Figure~\ref{fig:VMSscalability}.

\begin{figure}[t]
  \centering
  \includegraphics[bb=0 0 361 311,width=0.45\textwidth]{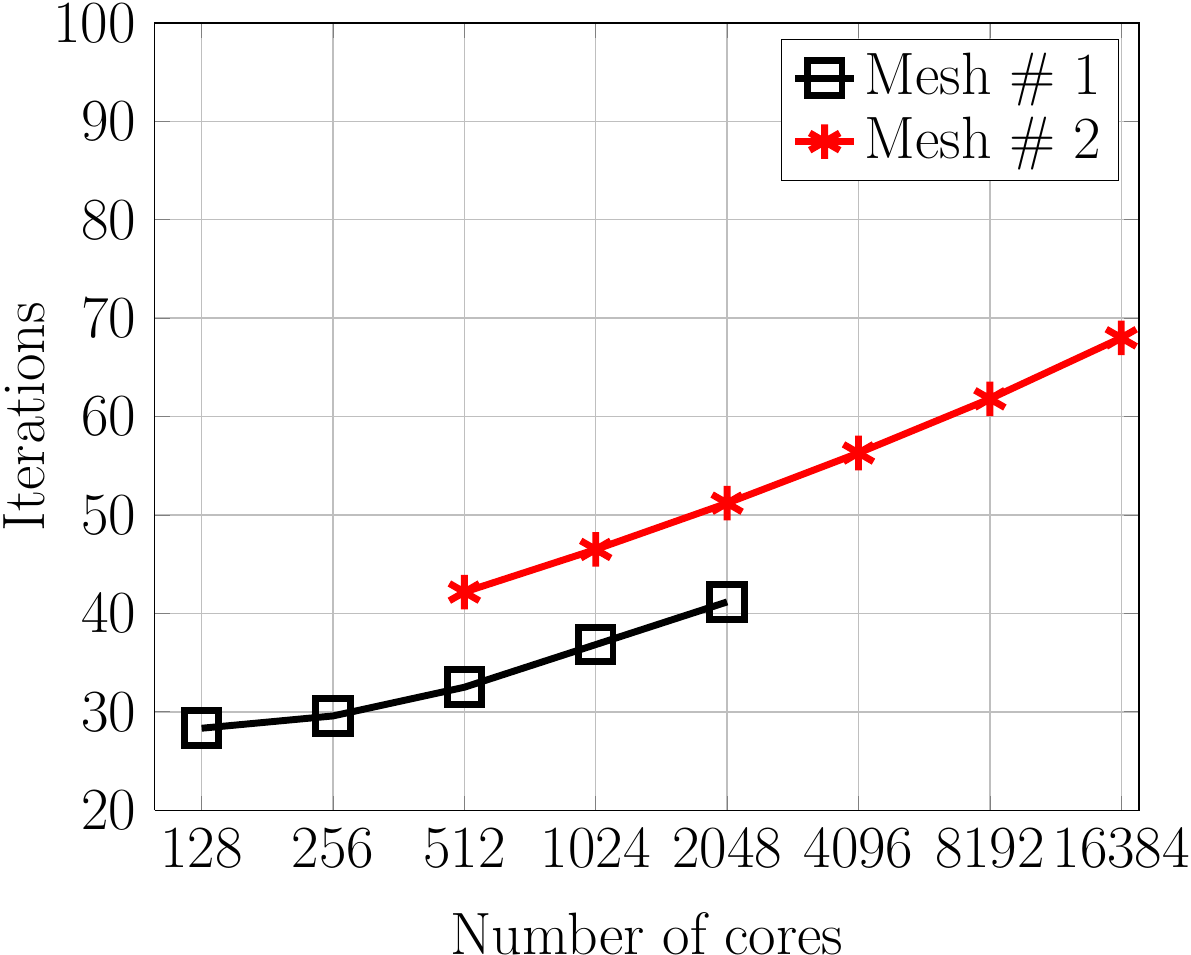}
  \includegraphics[bb=0 0 361 311,width=0.45\textwidth]{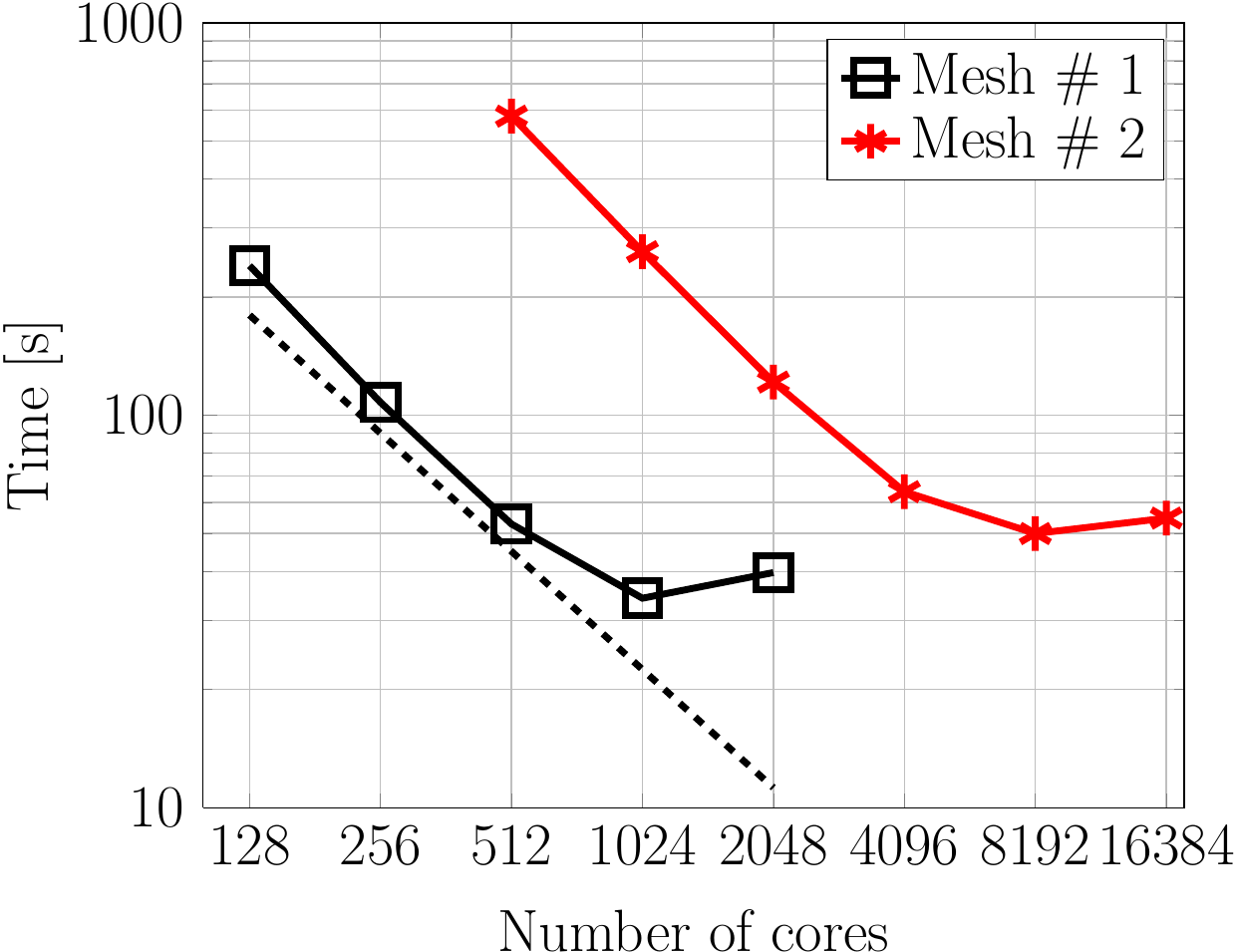}
  \caption{Blood flow in a patient-specific arterial bypass. Scalability of FaCSI preconditioner for FSI simulations~\cite{DeparisFortiGrandperrinQuarteroniFaCSI2015}.
  Simulations on PizDora at CSCS. Finite elements used: fluid P2-P1, structure P2, 
  ALE mapping P2. See Table~\ref{Tab:cylDofs} for the number of degrees of freedom.}
  \label{fig:FSIscalability}
\end{figure}

Applying the same techniques to build a parallel preconditioner for FSI 
based on a Dirichlet--Neumann splitting~\cite{CrosettoDeparisFouresteyQuarteroni2011}
and a SIMPLE~\cite{PatankarSpalding1972,Elman2008} preconditioner
does not lead to a scalable algorithm. To this end, it is necessary 
to add additional algebraic operations which leads to a FaCSI
preconditioner~\cite{DeparisFortiGrandperrinQuarteroniFaCSI2015}, see also Figure~\ref{fig:FSIscalability}.

\begin{table}%
\begin{tabular}{lccccc}
\hline
Mesh & Fluid DoF & Structure DoF & Coupling DoF & Geometry DoF & Total \\
\hline
Coarse & 9'029'128 & 2'376'720 & 338'295 & 8'674'950 & 20'419'093 \\
Fine & 71'480'725 & 9'481'350 & 1'352'020 & 68'711'934 & 151'026'029\\ 
\hline
\end{tabular}
\caption{Femoropopliteal bypass test case: number of Degrees of Freedom (DoF).}\label{Tab:cylDofs}
\end{table}

\subsection{Algebraic Factorizations}

Another issue is to reduce the computational cost by separating pressure and velocity computations. 
The origin of splitting schemes for Navier-Stokes problems may be dated back to the 
separate work of A. Chorin and R. Temam that give rise to the well-known scheme that bears their name.
The basic scheme operates at a differential level by exploiting the Helmholtz decomposition Theorem (also known as Ladhyzhenskaja Theorem)
to separate the differential problem into the sequence of
a vector advection-diffusion-reaction problem and a Poisson equation, with a final correction step for the velocity.
As opposed to the ``split-then-discretized'' paradigm, B. Perot in ~\cite{perot} advocates
a ``discretize then split'' strategy, by pointing out the formal analogy
between the Chorin-Temam scheme and an inexact LU factorization of the matrix obtained after 
discretization of the Navier-Stokes equations. 
This latter approach, often called ``algebraic factorization''
is easier to implement, particularly when one has to treat general boundary conditions \cite{QuarteroniSaleriVenezianiCMAME}. 

Let us introduce briefly a general framework. 
Let $A$ be the matrix obtained by the finite element discretization of the (linearized) incompressible Navier-Stokes equations.
The discretized problem at each step reads
$$
A 
\begin{bmatrix}
\U \\
\uu p
\end{bmatrix}= 
\begin{bmatrix}
\uu f \\
0
\end{bmatrix}
\qquad {\rm with \ } \quad
A = 
\begin{bmatrix}
C & D^T \\
D & 0
\end{bmatrix}
$$
where $A$ collects the contribution of the 
linearized differential operator acting on the velocity field in the momentum equation,
$D$ and $D^T$ are the discretization of the divergence  and the gradient operators, respectively. 
Notice that
$$
A = LU=
\begin{bmatrix}
C & 0 \\
D & - D C^{-1} D^T
\end{bmatrix}
\begin{bmatrix}
	I & C^{-1}D^T \\
	0 & I
\end{bmatrix}.
$$
This ``exact'' LU factorization of the problem formally realizes a velocity-pressure splitting.
However there is no computational advantage because of the presence of the matrix $C^{-1}$, which is
not explicitly available, so any matrix-vector product with this matrix requires to solve a linear system.
The basic idea of algebraic splitting  is to approximate this factorization. A first possibility is to replace $C^{-1}$ with the inverse of the velocity mass matrix scaled by $\Delta t$. This is the result of the first term truncation
of the Neumann expansion that may be used to represent $C^{-1}$.   The advantage of this approximation is that
the mass matrix can be further (and harmless) approximated by a diagonal matrix by the popular ``mass-lumping'' step. 
In this way, $D C^{-1} D^T$ is approximated by a s.p.d matrix --- sometimes called ``discrete Laplacian'' for its spectral analogy with the Laplace operator ---
that can be tackled with many different convenient numerical strategies.
In addition,  it is possible to see that the splitting error gathers in the
first block row, i.e.\ in the momentum equation. We finally note that replacing the original $C^{-1}$ with the velocity
mass matrix in the $U$ factor of the splitting implies that the exact boundary conditions cannot be enforced exactly.
The Yosida strategy, on the other hand, follows a similar pathway, except for not approximating 
$C^{-1}$ in $U$. Similar properties can be proved as for the Perot scheme, however in this case
the splitting error affects  only the mass conservation (with a moderate mass loss depending on the time step)
and the final step does actually enforce the exact boundary conditions for the velocity.
Successively, different splittings have been proposed in \cite{GauthierSaleriVeneziani,SaleriVeneziani,GervasioSaleriVeneziani,Gervasio2006,Gervasio2008,Veneziani2009}
to reduce the impact of the splitting error by  successive corrections of the pressure field. 
In particular, in \cite{Veneziani2003,GauthierSaleriVeneziani} the role of inexact factorizations as preconditioners
for the original problem  was investigated. 

LifeV incorporates these last developments. In particular, the Yosida
scheme has been preferred since the error on the mass conservation has
less impact on the interface with the structure in fluid-structure
interaction problems.  It is worth noting that a special block operator structure reflecting
the algebraic factorization concept has been implemented in
\cite{VillaPhD}.
 
It is worth noting that a robust validation of these methods has been
successfully performed not only against classical analytical test
cases but also within the framework of one {\it Critical Path
  Initiative} promoted by the US Food and Drug Administration (FDA)
\cite{QuainiEtAl}, (https://fdacfd.nci.nih.gov). Also, extensions of the inexact algebraic factorization
approach to the steady  problem have been recently proposed in \cite{Viguerie}.
 
\paragraph{Time adaptivity}

An interesting follow up of the pressure corrected Yosida algebraic factorizations is presented in  
\cite{VenezianiVilla2012}. This work stems from the fact that the sequence of pressure corrections
not only provides an enhancement of the overall  splitting error, but also provides an error estimator in time
for the pressure field - with no additional computational cost. Based on this idea,
a sophisticated time adaptive solver has been introduced~\cite{VillaPhD,VenezianiVilla2012},
with the aim of cutting  the computational costs by a smart and automatic selection of the time step.
The latter must be the trade-off among the desired accuracy, the computational efficiency and the
numerical stability constraints introduced by the splitting itself.
The final result is a solver that automatically detects the optimal time step, possibly performing
an appropriate number of pressure correction to attain stability.

This approach is particularly advantageous for computational hemodynamics problems featuring 
a periodic alternation of fast and slow transients (the so called systolic and diastolic phases in
circulation). As a matter of fact, for the same level of accuracy, the total number of time steps
required within a heart beat is reduced to one third of the ones required by the
non adaptive scheme.

In fact, in \cite{VenezianiVilla2012} a smart combination of algebraic factorizations
as solvers and preconditioners of the Navier-Stokes equations
based on the {\it a posteriori} error estimation provided by the pressure corrections
is proposed as a potential optimal trade-off between numerical stability and efficiency.


\section{Beyond the proof of concept}
\label{applications}

As explained in Sect. \ref{sec:introduction}, LifeV is intended to be a tool to work aggressively on real problems,
aiming at a general scope of bringing most advanced methods
for computational predictive tools in the engineering practice (in broad sense).

In particular, one of the most important applications --- yet non exclusive --- is the simulation of cardiovascular problems.
Examples of the use of LifeV for  real clinical problems are: simulations of Left Ventricular Assist Devices (LVAD)~\cite{bonnemain13:_numer,bonnemain12:VAD},
the study of the physiological~\cite{CrosettoReymondDeparisKontaxakisStergiopulosQuarteroni2011} and abnormal fluid-dynamics in ascending aorta
in presence of a bicuspid aortic valve~\cite{bonomiv1,vergarav1,faggianoa1,Viscardi2010},
of Thoracic EndoVascular Repair (TEVAR) \cite{pavia1,pavia2}, of the Total Cavopulmonary Connection
\cite{y1,restrepo2015surgical,tang2015respiratory}, of blood flow in stented coronary arteries \cite{jacc} and in cerebral aneurysms \cite{passerinietal}.
The current trend in this field is the setting up of {\it in silico} or ``Computer Aided'' Clinical Trials\cite{AV4Gogas},
i.e. of systematic investigations on large pools of patients to retrieve data of clinical relevance by integrating traditional measures and numerical simulations\footnote{The ABSORB Project granted by Abbott Inc. at Emory University and the iCardioCloud Project granted by Fondazione Cariplo to University of Pavia have been developed in this perspective using LifeV.}
In this scenario, numerical simulations are part of a complex, integrated pipeline involving: (i) Image/Data retrieval; (ii) Image Processing and Reconstruction (extracting the patient-specific morphology); (iii) Mesh generation and preprocessing (encoding of the boundary conditions);  
(iv) Numerical simulation (with LifeV); (v) Postprocessing and synthesis. As in the following sections we mainly describe applications relevant to step (iv), it is important to stress the integrated framework in which LifeV developers work, toward a systematic automation of the process, required from the large volume of patients to process. 
In this respect,
we may consider LifeV as a vehicle of {\it methodological transfer} or {\it translational mathematics}, as leading edge methods
are made available to the engineering community with a short time-to-market. Here
we present a series of distinctive applications where we feel that using LifeV actually allowed to
to bring rapidly new methods to real problems beyond the proof of concept stage.

\subsection{FSI}

In the Fluid-Structure Interaction (FSI) context, LifeV has offered a
very important bench for testing novel algorithms.
For example, it has been possible to test Robin-based interface conditions
for applications in hemodynamics~\cite{nobilev1,nobilep1,nobilep2},
or compare segregated algoritms, the monolithic formulation, and the Steklov-Poincar\'e formulation~\cite{DeparisDiscacciatiFouresteyQuarteroni2006}.

An efficient solution method for FSI problem considers the
physical unknowns and the fluid geometry problem for the Arbitrary
Lagrangian-Eulerian (ALE) mapping as a single variable.
This monolithic description implies to use ad-hoc parallel preconditioners.
Most often they rely on a Dirichlet-Neumann inexact factorization
between fluid and structure and then specific preconditioners
for the subproblems~\cite{CrosettoDeparisFouresteyQuarteroni2011}.
Recently a new preconditioner FaCSI~\cite{DeparisFortiGrandperrinQuarteroniFaCSI2015}
has been developed and tested
with LifeV with an effective scalability up to 4 thousands processors.
A next step in FSI has been the use of non-conforming
meshes between fluid and structure using rescaled-localized radial basis
functions
\cite{DeparisFortiQuarteroniRLRBF2014}.

A study on different material constitutive models for cerebral arterial tissue
--- in particular Hypereslatic isotropic laws, Hyperelastic anisotropic   laws ---
have been studied in~\cite{tricerri15:_fluid}.
A benchmark for the simulation of the flow inside carotids and
the computation of shear stresses \cite{FSIbenchmark2015} has been tested with
LifeV  coupled  with the FEAP library~\cite{feap},
which includes sophisticated anisotropic material models.

\subsection{Geometric Multiscale}

The cardiovascular system features coupled local and global dynamics.
Modeling its integrity by three dimensional geometries including FSI
is either unfeasible or very expensive computation-wise and,
most of the time, useless. A more efficient model entails
for the coupling of multiple dimensions, like lumped zero-dimensional
models, hyperbolic one-dimensional ones, and three dimensional FSI,
leading to the so called {\it geometric multiscale} modeling as advocated in
\cite{VenezianiPhD,formaggia2001coupling}. 

LifeV implements a full set of tools to integrate 0D, 1D and 3D-FSI models of the cardiovascular
system~\cite{malossi11:algorithms,malossi13:_implic,blanco:TotStress2013,passerini20093d},
with also a multirate time stepping scheme to improve the computational efficiency~\cite{malossi:two_time_step}.
It has been used to simulate integrated models of the cardiovascular system~\cite{bonnemain13:_numer}.

A critical aspect of this approach is the management of the dimensional mismatch
between the different models, as the accurate 3D problems
require more information at the interface that the one provided by the other
surrogate models, This required the accurate analysis of ``defective boundary problems''
~\cite{formaggia02,veneziani2005flow,venezianiv2,formaggia2008new,formaggia2010flow}. 
A recent review on these topics can be found in \cite{quarteroni2016geometric}.


\subsection{Heart dynamics}

Electrocardiology is one of the problems - beyond CFD but still related to cardiovascular mathematics -
 where LifeV has cumulated extensive experience.
An effective preconditioner for the bidomain equations  has been proposed and demonstrated in  \cite{giorda09}.
 The basic idea is to use the simplified extended monodomain model to precondition the solution
 of the more realistic bidomain equations. Successively, the idea has been adapted to reduce the computational costs
 by mixing Monodomain and Bidomain equations in an adaptive procedure. A suitabe
 {\it a posteriori} estimator is used to decide when the Monodoimain equations are enough or
 the bidomain solution is needed \cite{MirabellaNobileVeneziani,GGPeregoVeneziani1}. Ionic models
 solved in LifeV ranges from the classical Rogers McCulloch, Fenton Karma, Luo Rudy I and II~\cite{clayton2011models} to
 more involved ones \cite{dfgqr_mmas15}. Specific high order methods (extending the classical Rush Larsen one)
 have been proposed and implemented in the library \cite{PeregoVeneziani}.
 
In addition, one research line has been oriented to the coupling of electrocardiology with cardiac mechanics.
Hyperelasticity problems based on non-trivial mixed and primal formulations with applications in
cardiac biomechanics have been studied in \cite{rrpq_pamm11,rrpq_ijnmbe12}. LifeV-based simulations of fully coupled electromechanics (using modules for the abstract coupling of solvers)
can be found in \cite{nqr_ijnmbe12,raprq_springer13,rlrsq_ejmsol14,abqr_m3as15} for whole organ models, and in \cite{rgrclsq_mmb14,grrlcf_springer15,ruiz_jcp15} for single-cell problems. The coupling with ventricular fluid dynamics and arterial tree FSI description are possible through a multiscale framework \cite{quarte_rev15}.
 The coupling the Purkinje network, a network of high electrical conductivity myocardium fibers, has been implemented in~\cite{Vergara2016218}.

Another research line in this field successfully carried out with LifeV is the variational estimation
of cardiac conductivities (the tensor coefficients that are needed by the Bidomain equations)
from potential measures \cite{YangVeneziani}.

\subsection{Inverse problems and data assimilation}

One of the most recent challenging topics in computational hemodynamics is the
quantification of uncertainty and the improvement of the reliability in patient-specific settings.
As a matter of fact, while the inclusion of patient specific geometries
is now a well established procedure (as we recalled above), many other
aspects of the patient-specific modeling still deserve attention. Parameters like viscosity,
vessel wall rigidity, or cardiac conductivity are not routinely measured (or measurable)
in the specific patient and however have generally a major impact on the numerical results.
These concepts have been summarized in \cite{veneziani2013inverse}.
Variational procedures have been implemented in LifeV, where the assimilation
with available data or the parameter estimation are obtained by minimizing a mismatch functional.
In \cite{deliaperegoveneziani} this approach was introduced to incorporate into the numerical simulation
of the incompressible Navier-Stokes equations
sparse data available in the region of interest;
  in \cite{peregov1} the procedure was introduced for estimating the vascular rigidity
  by solving an inverse FSI problem, while a similar procedure in \cite{YangVeneziani,YangPOD}
  aims at the estimate of the cardiac conductivity.

\subsection{Model reduction}

One of the major challenges of modern scientific computing is the controlled reduction of the computational costs.
In fact, practical use of HPC demands for extreme efficiency --- even real time solutions. Improvement of computing architectures and
cloud solutions that make relatively easy the access to HPC facilties is only a partial answer to this need~\cite{guzzetti2016platform}.
From the modeling and methodological side,
we need also customized models that can realize the trade off between efficiency and accuracy.
These may be found by a smart combination of available High Fidelity solutions, according to the offline/online paradigm;
or by the inclusion of specific features of a problem that may bring a significant advantage in comparison with general versatile
but expensive methods. In LifeV these strategy have been both considered.
For instance, in \cite{colciago14:RFSI}, a model for blood flow dynamics in a fixed domain,
obtained by transpiration condition and a membrane model for the structure,
has been compared to a full three dimesional FSI simulation. The former model
is described only in the lumen with the Navier--Stokes equations, the
structure is taken into account
by surface Laplace--Beltrami operator on the surface representing the
fluid-strucutre interface. The reduced model allows for roughly one third
of the computational time and, in situations where the displacement of the
artery is pretty small, the dynamics, including e.g. wall shear stresses,
are very close if not indistinguishable from a full three dimensional simulation

In \cite{BertagnaVeneziani} a solution reduction procedure based on the Proper Orthogonal Decomposition
(POD) was used to accelerate the variational estimate of the Young modulus of vascular tissues
by solving an Inverse Fluid Structuire Interaction problem. POD wisely combines
available offline High Fidelity solutions to obtain a rapid (online) parameter estimation.
More challenging is the use of a similar approach for cardiac conductivities~\cite{yang2015parameter,YangPOD}, requiring nonstandard procedures.

A directional model reduction procedure called HiMod ({\it Hierarchical Model Reduction} - see
\cite{ernperottoveneziani,alettiperottoveneziani,perotto2014survey,perotto2014coupled,blanco2015hybrid,mansilla2016transversally}) to accelerate the computation of
advection diffusion reaction problems as well as incompressible fluids in pipe-like domains (or generally domains with a clear dominant direction,
like in arteries)
has been implemented in LifeV \cite{alettiperottoveneziani,guzzetti}.
These modules will be released in the library soon.

\subsection{Darcy equations and porous media}

Single and multi-phase flow simulators in fractured porous media are of
paramount importance in many fields like oil exploration and exploitation, CO$_2$
sequestration, nuclear waste disposal and geothermal reservoirs.  A single-phase
flow solver is implemented with the standard finite element spaces
Raviart-Thomas for the Darcy velocity and piecewise constant for the pressure. A
global pressure-total velocity formulation for the two-phase flow is developed
and presented in \cite{Fumagalli2011a,Fumagalli2012} where the equations are
solved using an IMPES-like technique. To handle fractures and faults in an
efficient and accurate way the extended finite element method is adopted to
locally enrich the cut elements, see \cite{Iori2011,DelPra2015a}, and a
one-codimensional problem for the flow is considered for these objects, see
\cite{Ferroni2014}.

\subsection{Ice sheets}
Another application that has used the LifeV library (at Sandia Nat Lab, Albuquerque, NM) is the simulation of  ice sheet flow. Ice behaves like a highly viscous shear-thinning incompressible fluid and can be modeled by nonlinear Stokes equations. In order to reduce computational costs several simplifications have been made to the Stokes model, exploiting the fact that ice sheets are very shallow. LifeV has been used to implement some of these models, including the Blatter-Pattyn (also known as First Order) approximation and the L1L2 approximation. The former model is a three-dimensional nonlinear elliptic PDE and the latter a depth-integrated integro-differential equations (see \cite{Perego2012}). In all models, nonlinearity has been solved with Newton method, coupling LifeV with Trilinos NOX package.
The LifeV based ice sheet implementation has been mentioned in \cite{Evans2012} as an example of modern solver design for the solution of earth system models. Further,
in \cite{Tezaur2015} the authors verified the results of another ice sheet code with those obtained using LifeV.

LifeV ice sheet module has been coupled with the climate library MPAS and used in inter-comparison studies to asses sensitivities to different boundary conditions  and forcing terms, see \cite{Edwards2014} and \cite{Shannon2013}. The latter (\cite{Shannon2013}) has been considered in the IPCC (Intergovernmental Panel on Climate Change) report of 2014.

In \cite{Perego2014} the authors perform a large scale PDE-constrained optimization to estimate the basal friction field (70K parameters) in Greenland ice sheet. For this purpose LifeV has been coupled with the Trilinos package ROL to perform  a reduced-gradient optimization using BFGS. The assembly of state and adjoint equations and the computation of objective functional and its gradient have been performed in LifeV.

\section{Perspectives}
LifeV has proved to be a versatile library for the study of numerical
techniques for large scale and multiphysics computations with finite
elements. The code is in continuous development. The latest
introduction of expression templates has increased the easiness of
usage at the high level. Unfortunately, due to the fact the way the code
has been developed, mainly by PhD students, not all the applications
have been ported yet to this framework. Work is ongoing in that
direction. The library has been coded using the C++98 
 yet porting is ongoing to exploit new features of the C++11 and C++14 standard, which can make the code more readable, user friendly
 and efficient.

As for the parallelization issues, LifeV relies strongly on the tools provided by the Trilinos libraries. We are following their development closely and
we will be ready to integrate all the new features the library will offer, in particular with respect to hybrid type parallelism.

A target use of the library the team is working on is running on cloud facilities \cite{slaw,guzzetti2016platform}
and GPL architectures.

\subsection*{Acknowledgments}
Besides the funding agencies cited in the introduction of this paper,
we have to acknowledge all the developers of LifeV, who have contributed
with new numerical methods, provided help to beginners, and struggled for
the definition of a common path. Porting and fixing bugs is also a remarkable task to
which they have constantly contributed.
It is difficult to name people,  the list would be anyway incomplete, a good representation
of the contributors is given in the text of the paper and, of course, on the developers website.
However, the three seniors authors of this paper wish to thank their collaborators
who over the years gave passion and dedication to the development.
SD and LF from Lausanne and Milan groups wish to thank
A. Quarteroni, who has contributed LifeV in many ways, by supporting its foundation,
by dedicating human resources to the common project, outreaching, and mentoring.
Also, worth of specific acknowledgment is G. Fourestey who has invested a lot of efforts to parallelize the code by introducing Trilinos
and many of the concepts described in the paper.
SD wishes to thank all his collaborators at CMCS in Lausanne, who have been able to work toghether in a very constructive environment. Particularly
P. Crosetto for introducing and testing the parallel framework for monolithic FSI and associated preconditioners;  C. Malossi for the design and implementation of the multiscale framework; S. Quinodoz for the design and implementation of the expression template module which is now at the core of the finite element formulation in LifeV; A. Gerbi, S. Rossi, R. Ruiz Baier, and P. Tricerri for their work on the mechanics of tissues, including active electromechainics;  G. Grandperrin and R. Popescu for thier studies on parallel algorithms
for the approximation of PDEs, in particolar the former and D. Forti for their
essential contribution for the Navier--Stokes equations and FSI;
 the promising and successfull coupling of LifeV with
 a reduced basis solver for 3D problems has been carried out by C. Colciago and N. Dal Santo.

LF wishes to thank specially Antonio Cervone, for his work on the phase field formulation for stratified fluids and the Stokes solver, Alessio Fumagalli for the Darcy solver, Guido Iori and Marco del Pra for the XFEM implementation for Darcy's flow in fractured media and the porting to C++11, Davide Baroli for his many contributions and the support of several Master students struggling with the code, Stefano Zonca for the DG-XFEM implementation of fluid structure interaction problems, Michele Lancellotti for his work on Large Eddy Simulation. And, finally, a thank to Christian Vergara who has exploited LifeV for several real life applications of simulations of the cardiovascular system, as well as giving an important contribution to the development of numerical techniques for defective boundary conditions and FSI.

AV wishes to thank all the Emory collaborators for their passion, dedication and courage
 joining him in the new challenging and enthusiastic adventure in Atlanta;  specifically
 (i) T. Passerini for the constant generous development
 and support not only in the development but also in the dissemination of LifeV among Emory students;
 (ii) M. Perego for his rigorous contribution and expansion of using
 LifeV at Sandia on the Ice-Sheet activity; (iii) U. Villa
 for the seminal and fundamental work in the Navier-Stokes modules that provided the basis of the current CFD in LifeV at Emory;
 (iv) L. Mirabella for bridging LifeV with the community of computational electrocardioologists and biomedical engineers at GA Tech.
The early  contribution of Daniele Di Pietro (University of Montpellier) for Expression Templates is gratefully acknowledged too.

Finally, all the authors together join to remember Fausto Saleri.
Fausto Saleri introduced the name "LiFE" many years ago, in a 2D, Fortran 77 serial code
for advection diffusion problems. Over the years, we made changes and additions,
yet we do hope to have followed his enthusiasm, passion and dedication to scientific computing.

\bibliographystyle{abbrvnat}
\bibliography{bibliography}




\end{document}